\documentclass[11pt]{amsart}

\usepackage{amscd}
\usepackage{amsmath}
\usepackage{graphicx}
\usepackage{amsfonts}
\usepackage{amssymb}
\begin{document}

\title[Lie Derivations]
{Characterization of  Lie Derivations on von Neumann Algebras}

\author{Xiaofei Qi}
\address
{Department of Mathematics\\
 Shanxi University\\
  Taiyuan 030006\\
   P. R.
of China} \email{qixf1980@126.com}
\author{Jinchuan Hou}
\address{Department of
Mathematics\\
Taiyuan University of Technology\\
 Taiyuan 030024\\
  P. R. of China;
Department of Mathematics\\
 Shanxi University\\
  Taiyuan 030006\\
   P. R.
of China } \email{jinchuanhou@yahoo.com.cn}

\thanks{This work is partially supported by National Natural Science
Foundation of China (11171249, 11101250) and Young Talents Plan for
Shanxi Univeristy.}

\begin{abstract}
Let ${\mathcal M}$ be a von Neumann algebra without central
summands of type $I_1$ and $\xi\in{\mathbb C}$ a scalar.  It is
shown that an additive map $L$ on $\mathcal M$ satisfies $L(AB-\xi
BA)=L(A)B-\xi BL(A)+L(B)A-\xi AL(B)$ whenever $A,B\in{\mathcal M}$
with $AB=0$ if and only if one of the following statements holds:
(1)  $\xi=1$,   $L=\varphi+f$, where $\varphi$ is an additive
derivation on $\mathcal M$ and $f$ is an additive map from
$\mathcal M$ into its center vanishing on $[A,B]$ with $AB=0$; (2)
  $\xi=0$,   $L(I)\in{\mathcal Z}({\mathcal M})$ and there
exists an additive derivation $\varphi$ such that
$L(A)=\varphi(A)+L(I)A$ for all $A$; (3) $\xi=-1$,  $L$ is a Jordan
derivation; (4)  $\xi$ is rational and $\xi\not=0, \pm1$, $L$ is an
additive derivation; (5) $\xi$ is not rational, there exists an
additive derivation $\varphi$ satisfying $\varphi(\xi I)=\xi L(I)$
such that $L(A)=\varphi(A) + L(I)A$ for all $A \in{\mathcal M}$. A
linear map $L$ on $\mathcal M$ satisfies $L(AB-\xi BA)=L(A)B-\xi
BL(A)+L(B)A-\xi AL(B)$ whenever $A,B\in{\mathcal M}$ with $AB=0$ if
and only if there exists a $T\in\mathcal M$ and a linear map
$f:{\mathcal M}\rightarrow{\mathcal Z}({\mathcal M})$ vanishing on
$[A,B]$ with $AB=0$ such that (i)   $\xi=1$, $L(A)=AT-TA+f(A)$ for
all $A\in\mathcal M$; (ii)  $\xi=0$, $L(I)\in{\mathcal Z}({\mathcal
M})$ and $L(A)=AT-(T-L(I))A$ for all $A\in{\mathcal M}$; (iii)
$\xi\not=0,1$, $L(A)=AT-TA$ for all $A \in{\mathcal M}$.

\end{abstract}
\maketitle

\section{Introduction}

Let $\mathcal R$ be an associative ring (or an algebra over a field
$\mathbb{F}$). Recall that an additive (a linear) map $\delta$ from
$\mathcal R$ into itself is called an additive (a linear) derivation
if $\delta(AB)=\delta(A)B+A\delta(B)$ for all $A$, $B\in {\mathcal
R}$. More generally, an additive (a linear) map $L$ from $\mathcal
R$ into itself is called an additive (a linear) Jordan derivation if
$L(AB+BA)=L(A)B+AL(B)+L(B)A+BL(A)$ for all $A,B \in \mathcal R$
(equivalently, $L(A^2)=L(A)A+AL(A)$ for all $A\in{\mathcal R}$ if
the characteristic of $\mathcal R$ is not 2); is called a Lie
derivation if $L([A,B])=[L(A),B]+[A,L(B)]$ for all $A,B \in \mathcal
R$, where $[A,B]=AB-BA$ is the Lie product of $A$ and $B$. The
problem of how to characterize the linear (additive) Jordan (Lie)
derivations of rings and algebras has received many mathematicians'
attention for many years.   Bre$\check{s}$ar in \cite{B} proved that
every additive Lie derivation on a prime ring $\mathcal R$ with
characteristic not 2 can be decomposed as $\tau + \zeta,$ where
$\tau$ is a derivation from $\mathcal R$ into its central closure
and $\zeta$ is an additive map of $\mathcal R$ into the extended
centroid $\mathcal{C}$ sending commutators to zero. Johnson \cite{J}
proved that every continuous linear Lie derivation from a
$C^*$-algebra $\mathcal A$ into a Banach $\mathcal A$-bimodule $M$
is standard, that is, can be decomposed as the form $\tau + h,$
where $\tau: \mathcal A \rightarrow M$ is a derivation and $h$ is a
linear map from $\mathcal A$ into the center of $M$ vanishing at
each commutator. Mathieu and Villena \cite{M} showed that every
linear Lie derivation on a $C^*$-algebra is standard. In \cite{QH}
Qi and Hou proved that the same is true for additive Lie derivations
of nest algebras on Banach spaces.  For other results, see
\cite{B2,H} and the references therein.

Recently, there have been a number of papers on the study of
conditions under which derivations of rings or operator algebras can
be completely determined by the action on some elements concerning
products (see \cite{B1,HQ1, WJ,L,P} and the references therein). For
Lie derivations, some works were also done.  A linear (an additive)
map $L: {\mathcal R}\rightarrow {\mathcal R}$ is said to be Lie
derivable at a point $Z$ if $L([A,B]) =[L(A),B] +[A,L(B)]$ for any
$A,B \in {\mathcal R}$ with $[A,B]=Z$. Clearly, this definition is
not valid for some $Z$, for instance, for $Z=I$, as the unit $I$ can
not be a commutator $[A,B]$ in general. It is also obvious that the
condition of maps Lie derivable at some point is much weaker than
the condition of being a Lie derivation. Qi and Hou \cite{QH3}
discussed such linear maps on ${\mathcal J}$-subspace lattice
algebras.  Lu and Jing in \cite{Lu} gave another kind of
characterization for Lie derivations as follows. Let $X$ be a Banach
space with $\dim X \geq 3$ and ${\mathcal B}(X)$ the algebra of all
bounded linear operators acting on $X$. It is proved in \cite{Lu}
that if $\delta:{\mathcal B}(X)\rightarrow{\mathcal B}(X)$ is a
linear map satisfying $\delta([A,B])=[\delta(A),B]+[A,\delta(B)]$
for any $A,B\in{\mathcal B}(X)$ with $AB=0$ (resp. $AB=P$, where $P$
is a fixed nontrivial idempotent), then $\delta=d+\tau$, where $d$
is a derivation of ${\mathcal B}(X)$ and $\tau:{\mathcal
B}(X)\rightarrow{\mathbb C}I$ is a linear map vanishing at
commutators $[A,B]$ with $AB=0$ (resp. $AB=P$). Later, this result
was generalized to the maps  on triangular algebras and prime rings
in \cite{JQ} and \cite{QH2} respectively. Since factor von Neumann
algebras are prime, as a consequence of the result for prime rings,
all additive maps $\delta$ on factor von Neumann algebras satisfying
$\delta([A,B])=[\delta(A),B]+[A,\delta(B)]$ for any $A,B$ with
$AB=0$ are characterized.  However the proof for factor von Neumann
algebras is not valid anymore for general von Neumann algebras. So,
it is natural to ask what happens when the concerned von Neumann
algebras are not factor.

Let $\mathcal A$ be an algebra over a field $\mathbb F$. For a
scalar $\xi\in{\mathbb F}$  and for $A,B\in{\mathcal A}$, if $AB=\xi
BA$, we say that $A$ commutes with $B$ up to a factor $\xi$. The
notion of commutativity up to a factor for pairs of operators is an
important concept and has been studied in the context of operator
algebras and quantum groups (ref. \cite{BBP, Ka}). Motivated by
this,  a binary operation $[A,B]_\xi =AB-\xi BA$, called $\xi$-Lie
product of $A$ and $B$, was introduced in  \cite{QH}.  An additive
(a linear) map $L:{\mathcal A} \rightarrow {\mathcal A}$ is called
an additive (a linear) $\xi$-Lie derivation if
$L([A,B]_{\xi})=[L(A),B]_{\xi}+[A,L(B)]_{\xi}$ for all $A,B \in
{\mathcal A}$. This conception unifies several well known notions.
It is clear that a $\xi$-Lie derivation is a derivation if $\xi =0$;
is a Lie derivation if $\xi=1$; is a Jordan derivation if $\xi=-1$.
The structure of $\xi$-Lie derivations was characterized in
triangular algebras and prime algebras in \cite{QH,QH4}
respectively. Particularly, we got a characterization of $\xi$-Lie
derivations on Banach space nest algebras and standard operator
algebras.

Thus, more generally, one may ask what is the structure of additive
(linear) maps $L$ that satisfy
$L([A,B]_\xi)=[L(A),B]_\xi+[A,L(B)]_\xi$ for any $A,B$ with $AB=0$?
The purpose of the present paper is to study this question for maps
on von Neumann algebras  and  characterize all such maps on general
von Neumann algebras. Note that every map on a commutative von
Neumann algebra is a Lie derivation. So it is reasonable to confine
our attention to the von Neumann algebras that have no central
summands of type $I_1$.

This paper is organized as follows. Let ${\mathcal M}$ be a von
Neumann algebra without central summands of type $I_1$. In Section
2, we deal with the case $\xi=1$, that is, the case of Lie product,
and show that every additive  map $L:{\mathcal
M}\rightarrow{\mathcal M}$ satisfies $L([A,B])=[L(A),B]+[A,L(B)]$
for any $A,B\in{\mathcal M}$ with $AB=0$ if and only if it has the
form $L=\varphi+f$, where $\varphi:{\mathcal M}\rightarrow{\mathcal
M}$ is an additive  derivation and $f:{\mathcal M}\rightarrow
{\mathcal Z}({\mathcal M})$, the center of $\mathcal M$, is an
additive  map vanishing on each commutator $[A,B]$ whenever $AB=0$
(Theorem 2.1).  Section 3 is devoted to discussing the case of
$\xi\neq 1$.  We show that every additive map $L:{\mathcal
M}\rightarrow{\mathcal M}$ satisfies that
$L([A,B]_\xi)=[L(A),B]_\xi+[A,L(B)]_\xi$ for any $A,B\in{\mathcal
M}$ with $AB=0$ if and only if  $L(I)\in{\mathcal Z}({\mathcal M})$
and,  (1) $\xi=0$, there exists an additive derivation $\varphi$
such that $L(A)=\varphi(A)+L(I)A$ for all $A\in{\mathcal M}$; (2)
$\xi=-1$, $L$ is a Jordan derivation, that is, $L(A^2)=L(A)A+AL(A)$
for all $A\in{\mathcal M}$; (3)  $\xi\in{\mathbb C}$ is rational and
$\xi\not=0,-1$, $L$ is an additive derivation; (4) $\xi\in{\mathbb
C}$ is not rational, there exists an additive derivation $\varphi$
satisfying $\varphi(\xi I)=\xi L(I)$ such that $L(A)=\varphi(A)
+L(I)A$ for all $A \in{\mathcal M}$ (Theorem 3.1). Moreover, in the
last case (4), if (i) $L$ is continuous when restricted on ${\mathbb
C}I$ or (ii) there exists a positive number $c$ and a subsequence of
integers $k_n\in {\mathbb Z}$ with $|k_n|\rightarrow\infty$ as
$n\rightarrow\infty$ such that $\|L(\xi^{k_n}I)\|\leq c|\xi|^{k_n}$,
then $L$ is an additive derivation (Corollary 3.3). Here we say that
a complex number is rational if it has rational real and imaginary
parts.  Particularly, we get a structure theorem for additive
$\xi$-Lie derivations on von Neumann algebras without central
summands of type $I_1$ (Corollary 2.5 and Corollary 3.4).

For the linear case, we show that, a linear map $L$ on $\mathcal M$
satisfies $L(AB-\xi BA)=L(A)B-\xi BL(A)+L(B)A-\xi AL(B)$ whenever
$A,B\in{\mathcal M}$ with $AB=0$ if and only if there exists a
$T\in\mathcal M$ and a linear map $f:{\mathcal
M}\rightarrow{\mathcal Z}({\mathcal M})$ vanishing on $[A,B]$ with
$AB=0$ such that (1) in the case $\xi=1$, $L(A)=AT-TA+f(A)$ for all
$A\in\mathcal M$; (2) in the case $\xi=0$, $L(I)\in{\mathcal
Z}({\mathcal M})$ and $L(A)=AT-(T-L(I))A$ for all $A\in{\mathcal
M}$; (3) in the case $\xi\not=0,1$, $L(A)=AT-TA$ for all $A
\in{\mathcal M}$ (See Theorem 2.1 and Corollary 3.2).

\section{Characterization of Lie derivations}

In this section, we consider the question of characterizing Lie
derivations by action at zero product on general von Neumann
algebras having no central summands of type $I_1$.

{\bf Theorem 2.1.} {\it Let ${\mathcal M}$ be a von Neumann algebra
without central summands of type $I_1$.   Suppose that $L:{\mathcal
M}\rightarrow{\mathcal M}$ is an additive   map. Then $L$ satisfies
that $L([A,B])=[L(A),B]+[A,L(B)]$ for any $A,B\in{\mathcal M}$ with
$AB=0$ if and only if there exists an additive  derivation
$\varphi:{\mathcal M}\rightarrow{\mathcal M}$ and an additive  map
$f:{\mathcal M}\rightarrow {\mathcal Z}({\mathcal M})$ that vanishes
each commutator $[A,B]$ whenever $AB=0$, such that
$L(A)=\varphi(A)+f(A)$ for all $A\in{\mathcal M}$, where ${\mathcal
Z}({\mathcal M})$ denotes the center of $\mathcal M$.

Particularly, $L$ is linear and satisfies that
$L([A,B])=[L(A),B]+[A,L(B)]$ for any $A,B\in{\mathcal M}$ with
$AB=0$ if and only if there exists some $T\in {\mathcal M}$ and a
linear map $f:{\mathcal M}\rightarrow {\mathcal Z}({\mathcal M})$
that vanishes each commutator $[A,B]$ whenever $AB=0$, such that
$L(A)=AT-TA+f(A)$ for all $A\in{\mathcal M}$.}

Before proving this theorem, we need some notations. We introduce
the concept of core-free projections, which had been defined in
\cite{Mi}. Let $\mathcal M$ be any von Neumann algebra and
$A\in{\mathcal M}$. Recall that the central carrier of $A$, denoted
by $\overline{A}$, is the intersection of all central projections
$P$ such that $PA=0$. If $A$ is self-adjoint, then the core of $A$,
denoted by $\underline{A}$, is sup$\{S\in{\mathcal Z}({\mathcal
M}):S=S^*, S\leq A\}$. Particularly, if $A=P$ is a projection, it is
clear that $\underline{P}$ is the largest central projection $\leq
P$. A projection $P$ is core-free if $\underline{P}=0$. It is easy
to see that $\underline{P}=0$ if and only if $\overline{I-P}=I$.

We first give two lemmas, which are needed in this paper.

{\bf Lemma 2.2.} (\cite{Mi}) {\it Let $\mathcal M$ be a von Neumann
algebra without central summands of type $I_1$. Then each nonzero
central projection $C\in{\mathcal M}$ is the carrier of a core-free
projection in $\mathcal M$. Particularly, there exists a nonzero
core-free projection $P\in{\mathcal M}$ with $\overline{P}=I$.}

We remark here that a little more can be said about the above lemma.
We in fact have that $\mathcal M$ is a von Neumann algebra without
central summands of type $I_1$ if and only if it has a projection
$P$ with $\underline{P}=0$ and $\overline{P}=I$.

\if If ${\mathcal M}={\mathcal A}\oplus{\mathcal M}_1$, where
$\mathcal A$ is the type $I_1$ summands of $\mathcal M$. Then
$\mathcal A$ is commutative. And $P=\left(\begin{array}{cc} P_0&0\\
0 & P_1\end{array}\right)$, where $P_0\in{\mathcal A}$ and
$P_1\in{\mathcal M}_1$. Since $P\geq\left(\begin{array}{cc} P_0&0\\
0 & 0\end{array}\right)\in{\mathcal A}\oplus 0\subseteq {\mathcal
Z}({mathcal M})$, we have $P_0=0$. So $P=\left(\begin{array}{cc}
0&0\\ 0 & P_1\end{array}\right)$. Thus $\overline{P}\leq
\left(\begin{array}{cc} 0&0\\ 0 & I_1\end{array}\right)<I$. A
contradiction. \fi

{\bf Lemma 2.3.} (\cite{Mi}) {\it Let $\mathcal M$ be a von Neumann
algebra. For projections $P,Q\in{\mathcal M}$, if
$\overline{P}=\overline{Q}\not=0$ and $P+Q=I$, then $T\in{\mathcal
M}$ commutes with $PXQ$ and $QXP$ for all $X\in{\mathcal M}$ implies
$T\in{\mathcal Z}({\mathcal M})$.}

{\bf Lemma 2.4.} {\it Let $\mathcal M$ be a von Neumann algebra.
Assume that $P\in{\mathcal M}$ is a projection with
$\underline{P}=0$ and $\overline{P}=I$. Then $P{\mathcal M}P\cap
{\mathcal Z}({\mathcal M})=(I-P){\mathcal M}(I-P)\cap {\mathcal
Z}({\mathcal M})=\{0\}$.}

{\bf Proof.} If $P\in{\mathcal M}$ is a projection such that
$\underline{P}=0$ and $\overline{P}=I$, then it is clear that
$\underline{I-P}=0$ and $\overline{I-P}=I$. So we need only to show
that $P{\mathcal M}P\cap {\mathcal Z}({\mathcal M})=\{0\}$. Assume
on the contrary that there is a nonzero element $A\in P{\mathcal
M}P\cap {\mathcal Z}({\mathcal M})$. Then there is a nonzero
projection $Q\in P{\mathcal M}P\cap {\mathcal Z}({\mathcal M})$. It
is clear that $0<Q\leq \underline{P}$, contradicting to the
assumption $\underline{P}=0$. \hfill$\Box$

{\bf Proof of Theorem 2.1.} By Lemma 2.2, we can find a non-central
core-free projection $P$ with central carrier $I$. In the sequel we
fix such a projection $P$. By the definitions of core and central
carrier, $I-P$ is core-free and $\overline{I-P}=I$. For the
convenience, denote ${\mathcal M}_{ij}=P_i{\mathcal M}P_j$,
$i,j\in\{1,2\}$, where $P_1=P$ and $P_2=I-P$. Then ${\mathcal
M}={\mathcal M}_{11}+{\mathcal M}_{12}+{\mathcal M}_{21}+{\mathcal
M}_{22}$. In all that follows, when  writing $S_{ij}$,  it always
indicates $S_{ij}\in {\mathcal M}_{ij}$.

The ``if'' part is obvious. We will prove the ``only if'' part by
checking several claims.

{\bf Claim 1.} $PL(P)P+(I-P)L(P)(I-P)\in{\mathcal Z}({\mathcal M})$.

For any $A_{12}\in {\mathcal M}_{12}$, since $A_{12}P=0$, we have
$L([A_{12},P])=[L(A_{12}),P]+[A_{12},L(P)]$. It follows that
$$L(A_{12})+L(A_{12})P-PL(A_{12})+A_{12}L(P)-L(P)A_{12}=0.$$ Writing
$L(A_{12})=S_{11}+S_{12}+S_{21}+S_{22}$,  we get
$$\begin{array}{rl}&S_{11}+A_{12}L(P)P+A_{12}(I-P)L(P)(I-P)\\&-PL(P)PA_{12}+2S_{21}+S_{22}-(I-P)L(P)PA_{12}=0.\end{array}$$
This implies that $A_{12}(I-P)L(P)(I-P)-PL(P)PA_{12}=0$. So
$$\begin{array}{rl}&PA(I-P)[(I-P)L(P)(I-P)+PL(P)P]\\
=&[(I-P)L(P)(I-P)+PL(P)P]PA(I-P)\end{array}\eqno(2.1)$$ for all
$A\in{\mathcal M}$.

Similarly, for any $(I-P)AP\in {\mathcal M}_{21}$, by using the
equation $(I-P)AP(I-P)=0$, one can show that
$$\begin{array}{rl}&(I-P)AP[(I-P)L(P)(I-P)+PL(P)P]\\
=&[(I-P)L(P)(I-P)+PL(P)P](I-P)AP\end{array}\eqno(2.2)$$ holds  for
all $A\in{\mathcal M}$.

By  Lemma 2.3,  Eq.(2.1) and Eq.(2.2) ensure that
$PL(P)P+(I-P)L(P)(I-P)\in{\mathcal Z}({\mathcal M})$, as desired.

 Denote $S=PL(P)(I-P)-(I-P)L(P)P$ and define a
map $\delta:{\mathcal M}\rightarrow{\mathcal M}$ by
$$\delta(A)=L(A)+SA-AS$$ for every $A\in{\mathcal M}$. It is easily checked that, for any $A,B\in{\mathcal
M}$,
$$AB=0\Rightarrow \delta([A,B])=[\delta(A),B]+[A,\delta(B)].$$
Moreover, by Claim 1, we have
$$\begin{array}{rl}\delta(P)=&L(P)-(I-P)L(P)P-PL(P)(I-P)\\
=&PL(P)P+(I-P)L(P)(I-P)\in{\mathcal Z}({\mathcal
M}).\end{array}\eqno(2.3)$$

{\bf Claim 2.} $P\delta(I)(I-P)=(I-P)\delta(I)P=0$.

Since $(I-P)P=0$, we have
$$0=\delta([I-P,P])=[\delta(I-P),P]+[I-P,\delta(P)]=[\delta(I),P],$$which
implies that $P\delta(I)(I-P)=(I-P)\delta(I)P=0$.

{\bf Claim 3.} $\delta({\mathcal M}_{ij})\subseteq{\mathcal
M}_{ij}$, $1\leq i\not=j\leq 2$.

We only need to check that $\delta({\mathcal
M}_{12})\subseteq{\mathcal M}_{12}$, and the proof of another
inclusion relation is similar.

For any $A_{12}\in {\mathcal M}_{12}$, write
$\delta(A_{12})=S_{11}+S_{12}+S_{21}+S_{22}$.  By the equation
$A_{12}P=0$, we have
$-\delta(A_{12})=\delta([A_{12},P])=[\delta(A_{12}),P]+[A_{12},\delta(P)]=[\delta(A_{12}),P]$.
It follows from Eq.(2.3) that $S_{11}+2S_{21}+S_{22}=0.$ This
implies that $S_{11}=S_{12}=S_{22}=0$, and so
$\delta(A_{12})=S_{12}\in{\mathcal M}_{12}$.

{\bf Claim 4.} There exists a map $f_i:{\mathcal
M}_{ii}\rightarrow{\mathcal Z}({\mathcal M})$   such that
$\delta(A_{ii})\in{\mathcal M}_{ii}+f_i(A_{ii})$ holds for any
$A_{ii}\in{\mathcal M}_{ii}$, $i=1,2$.

Take any $A_{11}\in {\mathcal M}_{11}$ and $B_{22}\in {\mathcal
M}_{22}$. Since $A_{11}B_{22}=0$, we have
$$0=\delta([A_{11},B_{22}])=[\delta(A_{11}),B_{22}]+[A_{11},\delta(B_{22})].$$
Let $\delta(A_{11})=S_{11}+S_{12}+S_{21}+S_{22}$ and
$\delta(B_{22})=T_{11}+T_{12}+T_{21}+T_{22}$. Then we get
$$(A_{11}T_{11}-T_{11}A_{11})+(S_{12}B_{22}+A_{11}T_{12})-(B_{22}S_{21}+T_{21}A_{11})+(S_{22}B_{22}-B_{22}S_{22})=0.$$
 It follows that
$$A_{11}T_{11}=T_{11}A_{11},\quad S_{22}B_{22}=B_{22}S_{22}\eqno(2.4)$$
and $$S_{12}B_{22}+A_{11}T_{12}=0,\quad
B_{22}S_{21}+T_{21}A_{11}=0.\eqno(2.5)$$ Fixing $A_{11}$ and letting
$B_{22}$ run over all ${\mathcal M}_{22}$, and in turn, fixing
$B_{22}$ and letting $A_{11}$ run over all ${\mathcal M}_{11}$, by
Eq.(2.4), we get
$$T_{11}=f_2(B_{22})P\ {\rm and}\  S_{22}=f_1(A_{11})(I-P)\ {\rm
hold \ for\ some}\ \ f_1(A_{11}), f_2(B_{22})\in{\mathcal
Z}({\mathcal M}).\eqno(2.6)$$ Taking $B_{22}=I-P$ in Eq.(2.5), then
for any $A_{11}$, by Eq.(2.3) and Claim 2, we have
$$S_{12}=-A_{11}P\delta(I-P)(I-P)=-A_{11}P\delta(I)(I-P)+A_{11}P\delta(P)(I-P)=0\eqno(2.7)$$
and
$$S_{21}=-(I-P)\delta(I-P)PA_{11}=-(I-P)\delta(I)PA_{11}+(I-P)\delta(P)PA_{11}=0.\eqno(2.8)$$
Similarly, taking $A_{11}=P$ in Eq.(2.5), one can easily check that
$T_{12}=T_{21}=0$. Combining this and Eqs.(2.6)-(2.8), we obtain
$$\delta(A_{11})=S_{11}+S_{22}=S_{11}- f_1(A_{11}) P+ f_1(A_{11})
I\in{\mathcal M}_{11}+ f_1(A_{11})$$ and
$$\delta(B_{22})=T_{11}+T_{22}=f_2(B_{22})
I-f_2(B_{22})(I-P)+T_{22}\in{\mathcal M}_{22}+f_2(B_{22}).$$  Hence
the Claim 4 is true.

Now let us define two maps  $f:{\mathcal M}\rightarrow{\mathcal
Z}({\mathcal M})$ and $d:{\mathcal M}\rightarrow{\mathcal M}$
respectively by
$$f(A)=f_1(A_{11})+f_2(A_{22})\ \mbox{and}\  d(A)=\delta(A)-f(A)$$ for
all $A=A_{11}+A_{12}+A_{21}+A_{22}\in{\mathcal M}.$ Then by Claims
3-4 we have
$$d({\mathcal M}_{ij})\subseteq{\mathcal M}_{ij},\
d({\mathcal M}_{ii})\subseteq{\mathcal M}_{ii},\ d({\mathcal
M}_{ij})=\delta({\mathcal M}_{ij}), \ 1\leq i\not= j\leq
2.\eqno(2.9)$$

{\bf Claim 5.} $d$ and $f$  are additive.

By the definition of $d$, $d$ is additive on ${\mathcal M}_{12}$ and
${\mathcal M}_{21}$. So we only need to verify that $d$ is additive
on ${\mathcal M}_{ii}$ $(i=1,2).$ In fact, for any $A_{11},B_{11}\in
{\mathcal M}_{11},$ we have $$\begin{array}{rl}
&d(A_{11}+B_{11})+f(A_{11}+B_{11})
=\delta(A_{11}+B_{11})=\delta(A_{11})+\delta(B_{11})\\
=&d(A_{11})+f(A_{11})+d(B_{11})+f(B_{11})
=d(A_{11})+d(B_{11})+(f(A_{11})+f(B_{11})),
\end{array}$$
that is,
$$d(A_{11}+B_{11})-(d(A_{11})+d(B_{11}))=(f(A_{11})
+f(B_{11})-f(A_{11}+B_{11}))\in{\mathcal Z}({\mathcal M}).$$ Since
${\mathcal M}_{11}\cap{\mathcal Z}({\mathcal M})=\{0\}$ by Lemma
2.4, we see that $d(A_{11}+B_{11})=d(A_{11})+d(B_{11})$, and
consequently, $f(A_{11})+f(B_{11})=f(A_{11}+B_{11})$.

Similarly, one can prove that $d$ is also additive on ${\mathcal
M}_{22}$.

{\bf Claim 6.} $d$ is a  derivation, that is, $d(AB)=d(A)B+Ad(B)$
for all $A,B\in{\mathcal M}$.

We will complete the proof of the claim by three steps.

{\bf Step 1.} $d(A_{ii}B_{ij})=d(A_{ii})B_{ij}+A_{ii}d(B_{ij})$ for
all $A_{ii}\in{\mathcal M}_{ii}$, $B_{ij}\in{\mathcal M}_{ij}$ and
$d(A_{ij}B_{jj})=d(A_{ij})B_{jj}+A_{ij}d(B_{jj})$ for all
$A_{ij}\in{\mathcal M}_{ij}$, $B_{jj}\in{\mathcal M}_{jj}$, $1\leq
i\not= j\leq 2$.

We only give the proof for the case $A_{11}\in{\mathcal M}_{11}$ and
$B_{12}\in{\mathcal M}_{12}$. The other cases can be dealt with
similarly.

In fact, for any $A_{11}\in{\mathcal M}_{11}$ and
$B_{12}\in{\mathcal M}_{12}$, since $B_{12}A_{11}=0$, by Eq.(2.9),
we have
$$\begin{array}{rl}-d(A_{11}B_{12})=&-\delta(A_{11}B_{12})=\delta([B_{12},A_{11}])\\
=&\delta(B_{12})A_{11}-A_{11}\delta(B_{12})+B_{12}\delta(A_{11})-\delta(A_{11})B_{12}\\
=&-A_{11}d(B_{12})-d(A_{11})B_{12},\end{array}$$ that is,
$d(A_{11}B_{12})=d(A_{11})B_{12}+A_{11}d(B_{12})$ for all
$A_{11}\in{\mathcal M}_{11}$ and $B_{12}\in{\mathcal M}_{12}$.

{\bf Step 2.} $d(A_{ii}B_{ii})=d(A_{ii})B_{ii}+A_{ii}d(B_{ii})$ for
all $A_{ii},B_{ii}\in{\mathcal M}_{ii}$, $i=1,2$.

Let $i\not=j$. For any $A_{ii},B_{ii}\in{\mathcal M}_{ii}$ and any
$S_{ij}\in{\mathcal M}_{ij}$, by Step 1, on the one hand, we have
$$d(A_{ii}B_{ii}S_{ij})=d(A_{ii}B_{ii})S_{ij}+A_{ii}B_{ii}d(S_{ij});$$on the other hand,
$$d(A_{ii}B_{ii}S_{ij})=d(A_{ii})B_{ii}S_{ij}+A_{ii}d(B_{ii}S_{ij})
=d(A_{ii})B_{ii}S_{ij}+A_{ii}d(B_{ij})S_{ij}+A_{ii}B_{ii}d(S_{ij}).$$
Comparing the above two equations, we see that
$$(d(A_{ii}B_{ii})-d(A_{ii})B_{ii}-A_{ii}d(B_{ii}))S_{ij}=0\eqno(2.10)$$
holds for all $S_{ij}\in{\mathcal M}_{ij}$.

Similarly, one can verify that
$$S_{ji}(d(A_{ii}B_{ii})-d(A_{ii})B_{ii}-A_{ii}d(B_{ii}))=0\eqno(2.11)$$
holds for all $S_{ji}\in{\mathcal M}_{ji}$. Also note that, by
Eq.(2.9), it is obvious that
$$S_{ij}(d(A_{ii}B_{ii})-d(A_{ii})B_{ii}-A_{ii}d(B_{ii}))=(d(A_{ii}B_{ii})-d(A_{ii})B_{ii}-A_{ii}d(B_{ii}))S_{ji}=0.$$
So it follows from Eqs.(2.10)-(2.11) and Lemma 2.3 that
$d(A_{ii}B_{ii})-d(A_{ii})B_{ii}-A_{ii}d(B_{ii})\in{\mathcal
Z}({\mathcal M})$, which implies that
$d(A_{ii}B_{ii})-d(A_{ii})B_{ii}-A_{ii}d(B_{ii})=0$ by Lemma 2.4.

Note that, by using Step 2 and the fact $d({\mathcal
M}_{ii})\subseteq {\mathcal M}_{ii}$ ($i=1,2$), one can get
$$d(P)=d(I-P)=0.\eqno(2.12)$$

{\bf Step 3.} $d(A_{ij}B_{ji})=d(A_{ij})B_{ji}+A_{ij}d(B_{ji})$ for
all $A_{ij}\in{\mathcal M}_{ij}$ and $B_{ji}\in{\mathcal M}_{ji}$,
$1\leq i\not= j\leq 2$.

Take any $A_{12}\in{\mathcal M}_{12}$ and $B_{21}\in{\mathcal
M}_{21}$. Since $(A_{12}B_{21}-A_{12}-B_{21}+(I-P))(P+B_{21})=0$, by
the definition of $d$, Claim 5 (additivity of $d$), Eq.(2.9) and
Eq.(2.12), we have
$$\begin{array}{rl}
&-d(A_{12}B_{21}-A_{12}+B_{21}A_{12}B_{21}-B_{21}A_{12})-f(A_{12}B_{21}-B_{21}A_{12})\\
=&-\delta(A_{12}B_{21}-A_{12}+B_{21}A_{12}B_{21}-B_{21}A_{12})\\
=&\delta([A_{12}B_{21}-A_{12}-B_{21}+(I-P),P+B_{21}])\\
=&[\delta(A_{12}B_{21}-A_{12}-B_{21}+I-P),P+B_{21}]
+[A_{12}B_{21}-A_{12}-B_{21}+I-P,\delta(P+B_{21})]\\
=&[d(A_{12}B_{21}-A_{12}-B_{21}+I-P),P+B_{21}]
+[A_{12}B_{21}-A_{12}-B_{21}+I-P,d(P+B_{21})]\\
=&-d(A_{12})B_{21}+d(A_{12})-B_{21}d(A_{12}B_{21})\\&+B_{21}d(A_{12})
-A_{12}d(B_{21})-d(B_{21})A_{12}B_{21}+d(B_{21})A_{12}.\end{array}$$
It follows from Step 1 that
$$\begin{array}{rl}
&d(A_{12}B_{21}-B_{21}A_{12})+f(A_{12}B_{21}-B_{21}A_{12})\\
=&d(A_{12})B_{21}-B_{21}d(A_{12})+A_{12}d(B_{21})-d(B_{21})A_{12}.\end{array}\eqno(2.13)$$
Multiplying by $A_{12}$ from the left side and the right side
respectively in Eq.(2.13), and applying  Eq.(2.9), we get
$$A_{12}d(B_{21}A_{12})-f(A_{12}B_{21}-B_{21}A_{12})A_{12}
=A_{12}B_{21}d(A_{12})+A_{12}d(B_{21})A_{12}$$ and
$$d(A_{12}B_{21})A_{12}+f(A_{12}B_{21}-B_{21}A_{12})A_{12}
=d(A_{12})B_{21}A_{12}+A_{12}d(B_{21})A_{12}.$$ These two equations,
together with Step 1, yield
$$-d(A_{12})B_{21}A_{12}-f(A_{12}B_{21}-B_{21}A_{12})A_{12}
=-d(A_{12}B_{21})A_{12}+A_{12}d(B_{21})A_{12}$$ and
$$-A_{12}B_{21}d(A_{12})+f(A_{12}B_{21}-B_{21}A_{12})A_{12}=
-A_{12}d(B_{21}A_{12})+A_{12}d(B_{21})A_{12}.$$ Comparing the above
two equations, one achieves
$$f(A_{12}B_{21}-B_{21}A_{12})A_{12}=0.\eqno(2.14)$$
Similarly, multiplying by $B_{21}$ from the left side and the right
side respectively in Eq.(2.13), one can verify
$$f(A_{12}B_{21}-B_{21}A_{12})B_{21}=0.\eqno(2.15)$$
Next we will prove $f(A_{12}B_{21}-B_{21}A_{12})=0$. To do this, for
any  $A_{12}$, let $A_{12}=V|A_{12}|$ be its polar decomposition.
Then Eq.(2.14) implies that
$f(A_{12}B_{21}-B_{21}A_{12})|A_{12}|=0$. So
$|A_{12}|f(A_{12}B_{21}-B_{21}A_{12})^*=0$. It follows that
$$A_{12}f(A_{12}B_{21}-B_{21}A_{12})^*=V|A_{12}|f(A_{12}B_{21}-B_{21}A_{12})^*=0.\eqno(2.16)$$
Similarly, one can show that
$$B_{21}f(A_{12}B_{21}-B_{21}A_{12})^*=0.\eqno(2.17)$$
Multiplying by $f(A_{12}B_{21}-B_{21}A_{12})^*$ in Eq.(2.13), by
using Eqs.(2.16)-(2.17), we get
$$\begin{array}{rl}&d(A_{12}B_{21}-B_{21}A_{12})f(A_{12}B_{21}-B_{21}A_{12})^*\\
&+f(A_{12}B_{21}-B_{21}A_{12})f(A_{12}B_{21}-B_{21}A_{12})^*=0.\end{array}\eqno(2.18)$$
Note that, by Eq.(2.9), Step 2 and Eqs.(2.16)-(2.17), we have
$$\begin{array}{rl}&d(A_{12}B_{21})f(A_{12}B_{21}-B_{21}A_{12})^*\\
=&d(A_{12}B_{21}Pf(A_{12}B_{21}-B_{21}A_{12})^*P)-A_{12}B_{21}d(Pf(A_{12}B_{21}-B_{21}A_{12})^*P)\\
=&-A_{12}B_{21}d(Pf(A_{12}B_{21}-B_{21}A_{12})^*P)\end{array}$$ and
$$\begin{array}{rl}&d(B_{21}A_{12})f(A_{12}B_{21}-B_{21}A_{12})^*\\
=&d(B_{21}A_{12}(I-P)f(A_{12}B_{21}-B_{21}A_{12})^*(I-P))\\&-B_{21}A_{12}d((I-P)f(A_{12}B_{21}-B_{21}A_{12})^*(I-P))\\
=&-B_{21}A_{12}d((I-P)f(A_{12}B_{21}-B_{21}A_{12})^*(I-P)).\end{array}$$
Hence Eq.(2.18) implies
$$f(A_{12}B_{21}-B_{21}A_{12})^*f(A_{12}B_{21}-B_{21}A_{12})f(A_{12}B_{21}-B_{21}A_{12})^*=0,$$
and so $f(A_{12}B_{21}-B_{21}A_{12})=0$. Thus Eq.(2.13) reduces to
$$d(A_{12}B_{21}-B_{21}A_{12})
=d(A_{12})B_{21}-B_{21}d(A_{12})+A_{12}d(B_{21})-d(B_{21})A_{12},$$
which implies that $d(A_{12}B_{21})=d(A_{12})B_{21}+A_{12}d(B_{21})$
and $d(B_{21}A_{12})=B_{21}d(A_{12})+d(B_{21})A_{12}$ hold for all
$A_{12}\in{\mathcal M}_{12}$ and $B_{21}\in{\mathcal M}_{21}$, as
desired.

Now by Steps 1-3, it is easily  checked that $d$ is an additive
derivation.

{\bf Claim 7.} $f([A,B])=0$ for all $A,B\in{\mathcal M}$ with
$AB=0$.

In fact, for any $A,B\in{\mathcal M}$ with $AB=0$, we have
$$\begin{array}{rl}
f([A,B])=&\delta([A,B])-d([A,B])\\
=&[\delta(A),B]+[A,\delta(B)]-d(AB-BA)\\
=&[d(A)-h(A),B]+[A,d(B)-h(B)]-d(AB-BA)\\
=&[d(A),B]+[A,d(B)]-d(AB-BA)=0.\end{array}$$

{\bf Claim 8.} The theorem holds.

Indeed, let $\varphi(A)=d(A)-(SA-AS)$ for all $A\in{\mathcal M}$;
then, by the definitions of $\delta$ and $d$, we have
$L(A)=\varphi(A)+f(A)$ for all $A\in{\mathcal M}$ . It is easy to
check that $\varphi$ is an additive derivation on $\mathcal M$.

Furthermore, if $L$ is linear, then $\varphi$ and $f$ are linear,
 too. As   linear derivations on von Neumann algebras are inner, we see that,
 there exists an element $T\in{\mathcal M}$ such that $\varphi(A)=AT-TA$ for all $A\in{\mathcal M}$.
 The proof is finished. \hfill$\Box$

By Theorem 2.1, we get a  characterization  of additive Lie
derivations immediately.

{\bf Corollary 2.5.} {\it Let ${\mathcal M}$ be a von Neumann
algebra without central summands of type $I_1$.   Suppose that
$L:{\mathcal M}\rightarrow{\mathcal M}$ is an additive map. Then $L$
is a Lie derivation if and only if there exists an additive
derivation $\varphi$ and an additive map $f:{\mathcal
M}\rightarrow{\mathcal Z}({\mathcal M})$ vanishing on each
commutator $[A,B]$ for any $A,B$ such that $L=\varphi+f$. }

\section{Characterization of $\xi$-Lie derivations}

In this section, we consider the question of characterizing
$\xi$-Lie derivations for $\xi\not=1$ by action at zero product on
von Neumann algebras.

{\bf Theorem 3.1.} {\it Let ${\mathcal M}$ be a von Neumann algebra
without central summands of type $I_1$.   Suppose that $L:{\mathcal
M}\rightarrow{\mathcal M}$ is an additive map and $\xi$ is a scalar
with $\xi\not=1$. Then $L$ satisfies that
$L([A,B]_\xi)=[L(A),B]_\xi+[A,L(B)]_\xi$ for any $A,B\in{\mathcal
M}$ with $AB=0$ if and only if $L(I)\in{\mathcal Z}({\mathcal M})$
and one of the following statements holds:}

(1) {\it $\xi\not=0,-1$,  there exists an additive derivation
$\varphi$ with $\varphi (\xi I)=\xi L(I)$ such that
$L(A)=\varphi(A)+L(I)A$  for all $A\in{\mathcal M}$; in
particular, $L$ is an additive derivation whenever $\xi$ is a
rational complex number.}

(2) {\it $\xi=0$, there exists an additive derivation $\varphi$ such
that   $L(A)=\varphi(A)+L(I)A$ for all $A\in{\mathcal M}$.}

(3) {\it  $\xi=-1$,  $L$ is an additive Jordan derivation, that is,
$L$ satisfies $L(A^2)=L(A)A+AL(A)$ for all $A\in{\mathcal M}$.}

For the linear maps, we have

{\bf Corollary 3.2.} {\it Let ${\mathcal M}$ be a von Neumann
algebra without central summands of type $I_1$.   Suppose that
$L:{\mathcal M}\rightarrow{\mathcal M}$ is a linear map and $\xi$ is
a scalar with $\xi\not=1$. Then $L$ satisfies that
$L([A,B]_\xi)=[L(A),B]_\xi+[A,L(B)]_\xi$ for any $A,B\in{\mathcal
M}$ with $AB=0$ if and only if one of the following statements
holds:}

(1) {\it $\xi\not=0$,  there is some $T\in \mathcal M$ such that
$L(A)=AT-TA$ for all $A \in{\mathcal M}$;}

(2) {\it $\xi=0$,  $L(I)\in{\mathcal Z}({\mathcal M})$ and there
exists some $T\in{\mathcal M}$ such that $L(A)=AT-(T-L(I))A$ for all
$A\in{\mathcal M}$.}


{\bf Proof.} Obviously, we need only check the ``only if'' part.
Note that every linear derivation on a von Neumann algebra is inner.
If $\xi\not=0,-1$, by Theorem 3.1, $L(I)=0$ and $L$ is a linear
derivation. Hence, there is an element $T\in{\mathcal M}$ such that
$L(A)=AT-TA$ for all $A$; if $\xi=-1$, by (3) of Theorem 3.1, $L$ is
a linear Jordan derivation. However every linear Jordan derivation
of a C$^*$-algebra is a derivation (ref. \cite[Teorem 2.4]{HL}).
Hence $L$ is a derivation and has the form $A\mapsto AT-TA$ for some
$T$. So the statement (1) of the corollary is valid. If $\xi=0$, by
Theorem 3.1, $L(I)\in{\mathcal Z}({\mathcal M})$ and there exists a
$T\in\mathcal M$ such that
 $L(A)=AT-TA+L(I)A=AT-(T-L(I))A$ for all
 $A$, that is, the statement (2) holds.  \hfill$\Box$

{\bf Proof of Theorem 3.1.} It is obvious that each of statements
(1), (2) and (3) implies that $AB=0\Rightarrow
L([A,B]_\xi)=[L(A),B]_\xi +[A,L(B)]_\xi$. For instance, assume that
(1) is valid. Then, for any $A,B$ with $AB=0$, we have
$$\begin{array}{rl} L([A,B]_\xi)=& -L(\xi BA)=-(\varphi(\xi BA)+\xi L(I)BA)\\
=&-(\varphi(\xi I)BA+\xi\varphi(BA)+ \xi
L(I)BA)=-(\xi\varphi(B)A+\xi B\varphi(A)+2\xi L(I)BA)\\
=& \varphi(A)B+A\varphi(B)+2L(I)AB-(\xi\varphi(B)A+\xi
B\varphi(A)+2\xi L(I)BA)\\
=&[L(A),B]_\xi +[A, L(B)]_\xi.
\end{array}$$

The following give a proof of the ``only if'' part. We use the same
symbols to that in the proof of Theorem 2.1. Particularly, $P$ is a
fixed projection in $\mathcal M$ with $\underline{P}=0$ and
$\overline{P}=I$. In the sequel, we always assume that $\xi\not=1$
and $L:{\mathcal M}\rightarrow{\mathcal M}$ is an additive map
satisfying $L([A,B]_\xi)=[L(A),B]_\xi+[A,L(B)]_\xi$ for
$A,B\in{\mathcal M}$ with $AB=0$. We will prove the ``only if" part
by several claims.

{\bf Claim 1.} $PL(I)(I-P)=(I-P)L(I)P=0$ and
$(I-P)L(P)(I-P)=PL(I-P)P=0$.

Since $P(I-P)=0$, we have $[L(P),I-P]_\xi+[P,L(I-P)]_\xi=0$, that
is,
$$L(P)(I-P)-\xi (I-P)L(P)+PL(I-P)-\xi L(I-P)P=0;\eqno(3.1)$$
since $(I-P)P=0$, we have
$[L(I-P),P]_\xi+[I-P,L(P)]_\xi=0$, that is,
$$L(I-P)P-\xi PL(I-P)+(I-P)L(P)-\xi L(P)(I-P)=0.\eqno(3.2)$$
Multiplying by $P$ and $I-P$ from the left and the right
respectively in Eq.(3.1), one gets $PL(P)(I-P)+PL(I-P)(I-P)=0$, and
so $$PL(I)(I-P)=0;$$ multiplying by $I-P$ and $P$ from the left and
the right respectively in Eq.(3.2), one gets
$(I-P)L(I-P)P+(I-P)L(P)P=0$, and so $$(I-P)L(I)P=0;$$ multiplying by
$I-P$ from both sides in Eq.(3.1), one gets
$(I-P)L(P)(I-P)-\xi(I-P)L(P)(I-P)=0$, and so $$(I-P)L(P)(I-P)=0;$$
multiplying by $P$ from both sides in Eq.(3.2), one gets
$PL(I-P)P-\xi PL(I-P)P=0$, and so $$PL(I-P)P=0.$$ Hence the claim is
true.

Now define a map $\delta:{\mathcal M}\rightarrow{\mathcal M}$ by
$\delta(A)=L(A)+SA-AS$ for each $A\in{\mathcal M}$, where
$S=PL(P)(I-P)-(I-P)L(P)P$. It is easy to verify that $\delta$ is an
additive map satisfying
$[\delta(A),B]_\xi+[A,\delta(B)]_\xi=\delta([A,B]_\xi)$ for
$A,B\in{\mathcal M}$ with $AB=0$. Moreover, by Claim 1, we also have
$$P\delta(I)(I-P)=(I-P)\delta(I)P=(I-P)\delta(P)(I-P)=P\delta(I-P)P=0.$$
Thus we get
$$\begin{array}{rl}\delta(P)=&L(P)+TP-PT=PL(P)P\\
=&P\delta(P)P-P(TP-PT)P=P\delta(P)P\end{array}\eqno(3.3)$$ and
$$\begin{array}{rl}\delta(I-P)=&L(I-P)+T(I-P)-(I-P)T=(I-P)L(I-P)(I-P)\\
=&(I-P)\delta(I-P)(I-P)-(I-P)(T(I-P)-(I-P)T)(I-P)\\=&(I-P)\delta(I-P)(I-P).\end{array}\eqno(3.4)$$

In the following, for the convenience,  we write $P_1=P$, $P_2=I-P$
and ${\mathcal M}_{ij}=P_i{\mathcal M}P_j$.

{\bf Claim 2.} $\delta({\mathcal M}_{ii})\subseteq {\mathcal
M}_{ii}$,   $i=1,2$.

For any $A_{11}\in{\mathcal M}_{11}$, since $A_{11}(I-P)=0$, we have
$[\delta(A_{11}),I-P]_\xi+[A_{11},\delta(I-P)]_\xi=0$. This and
Eq.(3.4) yield
$$\delta(A_{11})(I-P)-\xi(I-P)\delta(A_{11})=0.\eqno(3.5)$$ Multiplying
by $P$ from the left side in Eq.(3.5), one gets
$$P\delta(A_{11})(I-P)=0;\eqno(3.6)$$ multiplying by $I-P$ from both sides in
Eq.(3.5), one gets  $(1-\xi)(I-P)\delta(A_{11})(I-P)=0$, which
implies $$(I-P)\delta(A_{11})(I-P)=0\eqno(3.7)$$ as $\xi\not=1$. On
the other hand, since $(I-P)A_{11}=0$, we have
$[\delta(I-P),A_{11}]_\xi+[I-P,\delta(A_{11})]_\xi=0$, that is,
$(I-P)\delta(A_{11})-\xi\delta(A_{11})(I-P)=0$. Multiplying $P$ from
the right side in the equation, one gets
$$(I-P)\delta(A_{11})P=0.\eqno(3.8)$$
Combining Eqs.(3.6)-(3.8), we obtain $\delta(A_{11})\in{\mathcal
M}_{11}$. So $\delta({\mathcal M}_{11})\subseteq{\mathcal M}_{11}$.

The proof of $\delta({\mathcal M}_{22})\subseteq{\mathcal M}_{22}$
is similar and we omit it here.

{\bf Claim 3.} $\delta(I)\in{\mathcal Z}({\mathcal M})$ and
$\delta(P_i)A_{ij}=A_{ij}\delta(P_j)$ holds for any
$A_{ij}\in{\mathcal M}_{ij}$, $i\not=j\in\{1,2\}$.

Firstly, take any $A_{12}\in{\mathcal M}_{12}$.  Since
$A_{12}P_1=0$, we get
$$\begin{array}{rl}\delta(-\xi
A_{12})=&[\delta(A_{12}),P_1]_\xi+[A_{12},\delta(P_1)]_\xi\\
=&\delta(A_{12})P_1-\xi
P_1\delta(A_{12})+A_{12}\delta(P_1)-\xi\delta(P_1)A_{12};\end{array}\eqno(3.9)$$
since $P_2A_{12}=0$, we get
$$\begin{array}{rl}\delta(-\xi
A_{12})=&[\delta(P_2),A_{12}]_\xi+[P_2,\delta(A_{12})]_\xi\\
=&\delta(P_2)A_{12}-\xi
A_{12}\delta(P_2)+P_2\delta(A_{12})-\xi\delta(A_{12})P_2.\end{array}\eqno(3.10)$$
Eqs.(3.9) and (3.10) yield
$$\begin{array}{rl}&\delta(A_{12})P_1-\xi
P_1\delta(A_{12})+A_{12}\delta(P_1)-\xi\delta(P_1)A_{12}\\
=&\delta(P_2)A_{12}-\xi
A_{12}\delta(P_2)+P_2\delta(A_{12})-\xi\delta(A_{12})P_2.\end{array}\eqno(3.11)$$
If $\xi\not=0$, multiplying by $P_1$ and $P_2$ from the left and the
right respectively in Eq.(3.11), by Eqs.(3.3)-(3.4), one obtains
$$\delta(P_1)A_{12}=P_1\delta(P_1)P_1A_{12}=A_{12}P_2\delta(P_2)P_2=A_{12}\delta(P_2).$$
If $\xi=0$, by using of the relation
$(P_1+A_{12})(A_{12}-P_2)=(A_{12}-P_2)(P_1+A_{12})=0$, we have
$$(\delta(P_1)+\delta(A_{12}))(A_{12}-P_2)+(P_1+A_{12})(\delta(A_{12})-\delta(P_2))=0$$
and
$$(\delta(A_{12})-\delta(P_2))(P_1+A_{12})+(A_{12}-P_2)(\delta(P_1)+\delta(A_{12}))=0,$$
which mean that
$$\delta(P_1)A_{12}+\delta(A_{12})A_{12}-\delta(A_{12})P_2+P_1\delta(A_{12})+A_{12}\delta(A_{12})-A_{12}\delta(P_2)=0$$and
$$\delta(A_{12})P_1+\delta(A_{12})A_{12}+A_{12}\delta(A_{12})-P_2\delta(A_{12})=0.$$
Combining the above two equations, one obtains
$\delta(P_1)A_{12}=A_{12}\delta(P_2)$, which, together with
Eqs.(3.3)-(3.4), implies that
$P_1\delta(P_1)P_1A_{12}=A_{12}P_2\delta(P_2)P_2.$

Thus we have proved that
$$\delta(P_1)A_{12}=P_1\delta(P_1)P_1A_{12}=A_{12}P_2\delta(P_2)(P_2)=A_{12}\delta(P_2)\quad{\rm for \ \ all}\quad
A_{12}\in{\mathcal M}_{12}.$$ Note that
$P_1\delta(P_2)P_1=P_2\delta(P_1)P_2=0$. The above equation
implies$$P_1\delta(I)P_1A_{12}=A_{12}P_2\delta(I)P_2,$$and so
$$\begin{array}{rl}\delta(I)A_{12}=&(P_1\delta(I)P_1+P_2\delta(I)P_2)A_{12}\\
=&A_{12}(P_2\delta(I)P_2+P_1\delta(I)P_1)\\
=&A_{12}\delta(I)\end{array} \eqno(3.12)$$ holds for all
$A_{12}\in{\mathcal M}_{12}$.

Similarly, one can show that
$$\delta(P_2)A_{21}=A_{21}\delta(P_1)\ {\rm and}\ \delta(I)A_{21}=A_{21}\delta(I)\quad{\rm for \
\ all}\quad A_{21}\in{\mathcal M}_{21}.\eqno(3.13)$$ Now by Lemma
2.3, Eqs.(3.12) and (3.13), we get $\delta(I)\in{\mathcal
Z}({\mathcal M})$.

{\bf Claim 4.} For any $A_{ij}\in{\mathcal M}_{ij}$ ($1\leq
i\not=j\leq 2$), the following statements hold.

(1) If $\xi\not=-1$, then $\delta(A_{ij})\in{\mathcal M}_{ij}$.

(2) If $\xi=-1$, then $P_1\delta(A_{ij})P_1=P_2\delta(A_{ij})P_2=0$
and $\delta(A_{ij})A_{ij}+A_{ij}\delta(A_{ij})=0$.

To prove (1) we only need to check that $\delta(A_{12})\in{\mathcal
M}_{12}$ for all $A_{12}\in{\mathcal M}_{12}$, and the proof for
$A_{21}$ is similar.

For any $A_{12}\in{\mathcal M}_{12}$, Eq.(3.11) is true. Then,
multiplying by $P_1$ and $P_2$ from both sides respectively in
Eq.(3.11), and noting that $\xi\not=1$ and Eqs.(3.3)-(3.4), one can
easily check that
$$P_1\delta(A_{12})P_1=P_2\delta(A_{12})P_2=0.\eqno(3.14)$$

To complete  the proof of the statement (1), we have to check that
$P_2\delta(A_{12})P_1=0$. We will prove this by considering two
cases.

{\bf Case 1.} $\xi=0$.

For any $A_{12}$, since $A_{12}P_1=0$, we obtain
$$\delta(A_{12})P_1+A_{12}\delta(P_1)=0.\eqno(3.15)$$
Multiplying by $P_2$ from the left side in Eq.(3.15), one gets
$P_2\delta(A_{12})P_1=0$.  This and Eq.(3.14) yield
$\delta(A_{12})=P_1\delta(A_{12})P_2\in{\mathcal M}_{12}$.

{\bf Case 2.} $\xi\not=0,-1$.

In this case, take any $A_{12},B_{12}\in{\mathcal M}_{12}$. Since
$(B_{12}-P_2)(P_1+A_{12})=0$, by Eqs.(3.3)-(3.4), we have
$$\begin{array}{rl}&\delta(-\xi B_{12}+\xi A_{12})=\delta([B_{12}-P_2,P_1+A_{12}]_\xi)\\
=&[\delta(B_{12}-P_2),P_1+A_{12}]_\xi+[B_{12}-P_2,\delta(P_1+A_{12})]_\xi\\
=&\delta(B_{12})P_1+\delta(B_{12})A_{12}-\xi P_1\delta(B_{12})-\xi
A_{12}\delta(B_{12})+\xi A_{12}\delta(P_2)\\
&+B_{12}\delta(A_{12})-P_2\delta(A_{12})-\xi\delta(P_1)B_{12}-\xi\delta(A_{12})B_{12}+\xi\delta(A_{12})P_2.
\end{array}\eqno(3.16)$$
Multiplying by $P_2$ from both sides in Eq.(3.16) and applying
Eq.(3.14), we get
$$P_2\delta(B_{12})P_1A_{12}=\xi P_2\delta(A_{12})P_1B_{12}\quad{\rm for \ \ all}\quad A_{12},B_{12}\in{\mathcal M}_{12}.\eqno(3.17)$$
Note that, multiplying by $P_2$ and $P_1$ from the left and the
right respectively in Eq.(3.9), one sees that
$$P_2\delta(A_{12})P_1=P_2\delta(-\xi A_{12})P_1\quad{\rm for \ \ all}\quad A_{12}\in{\mathcal M}_{12}.\eqno(3.18)$$
It follows from Eqs.(3.17) and (3.18) that
$$\begin{array}{rl}&-P_2\delta(B_{12})P_1A_{12}=P_2\delta(\xi B_{12})P_1A_{12}\\
=&\xi P_2\delta(A_{12})P_1(\xi B_{12}) =\xi^2
P_2\delta(A_{12})P_1B_{12}\end{array}$$ for  all
$A_{12},B_{12}\in{\mathcal M}_{12}$. This and Eq.(3.17) imply that
$$P_2\delta(A_{12})P_1B_{12}=0\quad{\rm for \ \ all}\quad A_{12},B_{12}\in{\mathcal M}_{12}.\eqno(3.19)$$

Similarly,  multiplying by $P_1$ from both sides in Eq.(3.16), by
using Eqs.(3.14) and (3.18), one can show that
$$B_{12}P_2\delta(A_{12})P_1=0\quad{\rm for \ \ all}\quad A_{12},B_{12}\in{\mathcal M}_{12}.\eqno(3.20)$$
Also note that
$P_2\delta(A_{12})P_1B_{21}=B_{21}P_2\delta(A_{12})P_1=0$ for all
$B_{21}\in{\mathcal M}_{21}$. Then it follows from Lemma 2.3 and
Eqs.(3.19)-(3.20) that $P_2\delta(A_{12})P_1\in{\mathcal
Z}({\mathcal M})$, and hence $P_2\delta(A_{12})P_1=0$. So the
statement (1) holds.

To prove the statement (2), note  that, as $\xi=-1$, $\delta$ in
fact satisfies
$$AB=0\Rightarrow \delta(AB+BA)=\delta(A)B+A\delta(B)+\delta(B)A+B\delta(A).$$
Then, for any $A_{12}\in {\mathcal M}_{12}$, since
$(P_1+A_{12})(A_{12}-P_2)=0$, we have
$$\begin{array}{rl}&\delta(P_1+A_{12})(A_{12}-P_2)+(A_{12}-P_2)\delta(P_1+A_{12})\\
&+\delta(A_{12}-P_2)(P_1+A_{12})+(P_1+A_{12})\delta(A_{12}-P_2)=0.\end{array}$$
It follows from Claims 2-3 and Eq.(3.11) that
$\delta(A_{12})A_{12}+A_{12}\delta(A_{12})=0$, which, together with
Eq.(3.14), imply that the statement (2) is true.

{\bf Claim  5.} The following statements hold.

(1) If $\xi\not=0,-1$, then $\delta(\xi AB)=\xi \delta(A)B+\xi
A\delta(B)$ for all $A,B\in{\mathcal M}$.

(2) If $\xi=0$, then there exists an additive derivation $\varphi$
such that $\delta(A)=\varphi(A) +\delta(I)A$ for all $A\in{\mathcal
M}$.

(3) If $\xi=-1$, then $\delta(A^2)=\delta(A)A+A\delta(A)$ for all
$A\in{\mathcal M}$, that is, $\delta$ is an additive Jordan
derivation.

We will prove the claim by considering three cases.

{\bf Case 1.} $\xi\not=0,-1$.

In this case, we will show that $\delta(\xi AB)=\xi \delta(A)B+\xi
A\delta(B)$ for all $A,B\in{\mathcal M}$ by three steps.

{\bf Step 1.} $\delta(\xi A_{ii}B_{ij})=\xi\delta(A_{ii})B_{ij}+\xi
A_{ii}\delta(B_{ij})$ for all $A_{ii}\in{\mathcal M}_{ii}$,
$B_{ij}\in{\mathcal M}_{ij}$, $1\leq i\not= j\leq 2$.

In fact, for any $A_{ii}\in{\mathcal M}_{ii}$ and
$B_{ij}\in{\mathcal M}_{ij}$, since $B_{ij}A_{ii}=0$, by Claims 2
and (1)of Claim 4, we have
$$\begin{array}{rl}&-\delta(\xi A_{ii}B_{ij})=\delta([B_{ij},A_{ii}]_\xi)\\
=&\delta(B_{ij})A_{ii}-\xi A_{ii}\delta(B_{ij})+B_{ij}\delta(A_{ii})-\xi\delta(A_{ii})B_{ij}\\
=&-\xi A_{ii}\delta(B_{ij})-\xi\delta(A_{ii})B_{ij},\end{array}$$
that is, $\delta(\xi A_{ii}B_{ij})=\xi\delta(A_{ii})B_{ij}+\xi
A_{ii}\delta(B_{ij})$ for all $A_{ii}\in{\mathcal M}_{ij}$ and
$B_{ij}\in{\mathcal M}_{ij}$.

Similarly, one can check  the following.

{\bf Step 2.} $\delta(\xi A_{ij}B_{jj})=\xi\delta(A_{ij})B_{jj}+\xi
A_{ij}\delta(B_{jj})$ for all $A_{ij}\in{\mathcal M}_{ij}$,
$B_{jj}\in{\mathcal M}_{jj}$, $1\leq i\not= j\leq 2$.

{\bf Step 3.} $\delta(\xi A_{ii}B_{ii})=\xi\delta(A_{ii})B_{ii}+\xi
A_{ii}\delta(B_{ii})$ for all $A_{ii},B_{ii}\in{\mathcal M}_{ii}$,
$i=1,2$.

Let $i\not=j$. For any $A_{ii},B_{ii}\in{\mathcal M}_{ii}$ and any
$S_{ij}\in{\mathcal M}_{ij}$, by Step 1, on the one hand, we have
$$\delta(\xi A_{ii}B_{ii}S_{ij})=\xi\delta(A_{ii}B_{ii})S_{ij}+\xi A_{ii}B_{ii}\delta(S_{ij});$$
on the other hand,
$$\begin{array}{rl}\delta(\xi A_{ii}B_{ii}S_{ij})=&\xi\delta(A_{ii})B_{ii}S_{ij}+\xi
A_{ii}\delta(B_{ii}S_{ij})\\
=&\xi\delta(A_{ii})B_{ii}S_{ij}+\xi^2
A_{ii}\delta(\xi^{-1}B_{ii})S_{ij}+\xi
A_{ii}B_{ii}\delta(S_{ij}).\end{array}$$ Comparing the above two
equations, we see that
$$(\delta(A_{ii}B_{ii})-\delta(A_{ii})B_{ii}-\xi A_{ii}\delta(\xi^{-1}B_{ii}))S_{ij}=0,$$
which implies that
$$(\delta(\xi A_{ii}B_{ii})-\xi\delta(A_{ii})B_{ii}-\xi A_{ii}\delta(B_{ii}))S_{ij}=0\eqno(3.21)$$
holds for all $S_{ij}\in{\mathcal M}_{ij}$.

Similarly, one can verify that
$$S_{ji}(\delta(\xi A_{ii}B_{ii})-\xi\delta(A_{ii})B_{ii}-\xi A_{ii}\delta(B_{ii}))=0\eqno(3.22)$$
holds for all $S_{ji}\in{\mathcal M}_{ji}$. Also note that, by Claim
2, it is obvious that
$$S_{ij}(\delta(\xi A_{ii}B_{ii})-\xi\delta(A_{ii})B_{ii}-\xi A_{ii}\delta(B_{ii}))
=(\delta(\xi A_{ii}B_{ii})-\xi\delta(A_{ii})B_{ii}-\xi
A_{ii}\delta(B_{ii}))S_{ji}=0.$$ So it follows from Lemma 2.3 and
Eqs.(3.21)-(3.22) that $\delta(\xi
A_{ii}B_{ii})-\xi\delta(A_{ii})B_{ii}-\xi
A_{ii}\delta(B_{ii})\in{\mathcal Z}({\mathcal M})$, which implies,
by Lemma 2.4, that $\delta(\xi
A_{ii}B_{ii})-\xi\delta(A_{ii})B_{ii}-\xi A_{ii}\delta(B_{ii})=0$.

{\bf Step 4.} $\delta(\xi A_{ij}B_{ji})=\xi\delta(A_{ij})B_{ji}+\xi
A_{ij}\delta(B_{ji})$ for all $A_{ij}\in{\mathcal M}_{ij}$ and
$B_{ji}\in{\mathcal M}_{ji}$, $1\leq i\not= j\leq 2$.

For any $A_{ij}\in{\mathcal M}_{ij}$ and $B_{ji}\in{\mathcal
M}_{ji}$ with $i\not=j$, since
$(A_{ij}B_{ji}-A_{ij}-B_{ji}+P_j)(P_i+B_{ji})=0$, by the definition
of $\delta$, we have
$$-\delta(\xi A_{ij}B_{ji}-\xi A_{ij}+\xi B_{ji}A_{ij}B_{ji}-\xi B_{ji}A_{ij})
=\delta([A_{ij}B_{ji}-A_{ij}-B_{ji}+P_j,P_i+B_{ji}]_\xi).$$ Thus by
Claim 2 and (1) of Claim 4, the above equation   reduces to
$$\begin{array}{rl}
&\delta(\xi A_{ij}B_{ji})-\delta(\xi A_{ij})-\delta(\xi
B_{ji}A_{ij})\\
=&\delta(A_{ij})B_{ji}+\delta(B_{ji})P_i-\delta(P_j)B_{ji}-\xi\delta(A_{ij})\\
&-\xi
B_{ji}\delta(A_{ij})-A_{ij}B_{ji}\delta(P_i)+A_{ij}\delta(B_{ji})\\
&+B_{ji}\delta(P_i)-P_j\delta(B_{ji})-\xi\delta(P_i)A_{ij}-\xi
\delta(B_{ji})A_{ij}.\end{array}$$ Multiplying by $P_j$ from both
sides in the above equation, by Claims 2 and 4 again, one obtains
$\delta(\xi B_{ji}A_{ij})=\xi\delta(B_{ji})A_{ij}+\xi
B_{ji}\delta(A_{ij})$, as desired.

Now, for any $A,B\in{\mathcal M}$, by Steps 1-4 and the additivity
of $\delta$,  it is easily  checked that $\delta(\xi
AB)=\xi\delta(A)B+\xi A\delta(B)$ holds  for all $A,B\in{\mathcal
M}$. So the statement (1) of Claim 5 is true.

{\bf Case 2.} $\xi=0$.

In this case, $\delta$ satisfies $AB=0\Rightarrow
\delta(A)B+A\delta(B)=0$. We first show that
$$\delta(AB)=\delta(A)B+A\delta(B)-\delta(I)AB$$ for all
$A,B\in{\mathcal M}$.

Let $1\leq i\not= j\leq 2$. By Claim 2 and  (1) in Claim 4, the
relation $(A_{ii}+A_{ii}B_{ij})(P_j-B_{ij})=0$ entails that
$$\delta(A_{ii}B_{ij})=\delta(A_{ii})B_{ij}+A_{ii}\delta(B_{ij})-A_{ii}B_{ij}\delta(P_j)\eqno(3.23)$$
holds for any $A_{ii}\in{\mathcal M}_{ii}$ and $B_{ij}\in{\mathcal
M}_{ij}$;  the relation $(P_i-A_{ij})(B_{jj}+A_{ij}B_{jj})=0$
implies that
$$\delta(A_{ij}B_{jj})=\delta(A_{ij})B_{jj}+
A_{ij}\delta(B_{jj})-\delta(P_i)A_{ij}B_{jj}\eqno(3.24)$$ holds for
any $A_{ij}\in{\mathcal M}_{ij}$ and $B_{jj}\in{\mathcal M}_{jj}$.
Then, by Claim 3, Eq.(3.23) and using a similar argument to that of
Step 3 in Case 1, one can show that
$$\delta(A_{ii}B_{ii})=\delta(A_{ii})B_{ii}+A_{ii}\delta(B_{ii})-A_{ii}B_{ii}\delta(P_i)\eqno(3.25)$$
holds for any $A_{ii},B_{ii}\in{\mathcal M}_{ii}$. Next,  by the
equation $(A_{ij}+A_{ij}B_{ji})(P_i-B_{ji})=0$ and Claim 2 and (1)
in Claim 4, one can obtain that
$$\delta(A_{ij}B_{ji})=\delta(A_{ij})B_{ji}+A_{ij}\delta(B_{ji})-A_{ij}B_{ji}\delta(P_i)\eqno(3.26)$$
holds for any $A_{ij}\in{\mathcal M}_{ij}$ and $B_{ji}\in{\mathcal
M}_{ji}$. Finally, the additivity of $\delta$, together with
Eqs.(3.23)-(3.26), ensures that
$\delta(AB)=\delta(A)B+A\delta(B)-\delta(I)AB$ holds for any $A,B\in
{\mathcal M}$.

Now let $\varphi :{\mathcal M}\rightarrow{\mathcal M}$ be the map
defined by $\varphi(A)=\delta(A)-\delta(I)A$. Note that
$\delta(I)\in{\mathcal Z}({\mathcal M})$. Thus we have that
$$\begin{array}{rl}
\varphi(AB)=&\delta(AB)-\delta(I)AB  =
\delta(A)B+A\delta(B)-2\delta(I)AB \\
=&
(\varphi(A)+\delta(I)A)B+A(\varphi(B)+\delta(I)B)-2\delta(I)AB=\varphi(A)B+A\varphi(B)
\end{array}
$$
holds for any $A,B\in\mathcal M$. So $\varphi$ is an additive
derivation, and $\delta(A)=\varphi(A)+\delta(I)A$ for all $A$.
particularly, $\delta$ is a generalized derivation.

{\bf Case 3.} $\xi=-1$.

In this case, $\delta$ satisfies
$$AB=0\Rightarrow \delta(BA)=\delta(A)B+A\delta(B)+\delta(B)A+B\delta(A).$$
We will show that $\delta$ is a Jordan derivation, and therefore the
statement (3) holds.

Let $1\leq i\not= j\leq 2$. For any $A_{ii}\in{\mathcal M}_{ii}$ and
$B_{ij}\in{\mathcal M}_{ij}$, since $B_{ij}A_{ii}=0$, by Claim 2 and
(2) in Claim 4, one can verify
$$\delta(A_{ii}B_{ij})=\delta(A_{ii})B_{ij}+A_{ii}\delta(B_{ij})+\delta(B_{ij})A_{ii};\eqno(3.27)$$
for any $A_{ij}\in{\mathcal M}_{ij}$ and $B_{jj}\in{\mathcal
M}_{jj}$, by using of the relation $B_{jj}A_{ij}=0$, Claim 2 and (2)
in Claim 4, one can verify
$$\delta(A_{ij}B_{jj})=\delta(A_{ij})B_{jj}+
A_{ij}\delta(B_{jj})+B_{jj}\delta(A_{ij}).\eqno(3.28)$$ For any
$A_{ii},B_{ii}\in{\mathcal M}_{ii}$, by Claim 3, Eq.(3.27) and using
a similar argument to that of Step 3 in Case 1, one can show that
$$\delta(A_{ii}B_{ii})=\delta(A_{ii})B_{ii}+A_{ii}\delta(B_{ii}).\eqno(3.29)$$
For any $A_{ij}\in{\mathcal M}_{ij}$ and $A_{ji}\in{\mathcal
M}_{ji}$, since
$(A_{ij}A_{ji}+A_{ij}+A_{ji}+P_j)(P_i-A_{ij}-A_{ji}+A_{ji}A_{ij})=0$,
by Claim 2 and (2) of Claim 4, it is easily checked that
$$\delta(A_{ij}A_{ji})=\delta(A_{ij})A_{ji}+A_{ij}\delta(A_{ji})\ {\rm and}\
\delta(A_{ji}A_{ij})=\delta(A_{ji})A_{ij}+A_{ji}\delta(A_{ij})\eqno(3.30)$$
Now, combining Eqs.(3.27)-(3.30), it is easy to verify that
$\delta(A^2)=\delta(A)A+A\delta(A)$ for all $A\in{\mathcal M}$, that
is, $\delta$ is a Jordan derivation.

{\bf Claim 6.} If $\xi\neq 0,-1$, then there exists an additive
derivation $\varphi$ satisfying $\varphi(\xi I)=\xi\delta(I)$ such
that   $\delta( A)= \varphi(A)+ \delta(I)A$ holds for any $A$; in
particular, in the case that $\xi$ is a rational complex number,
$\delta$ is an additive derivation.

By (1) of Claim 5, $\delta(\xi AB)=\xi(\delta(A)B+A\delta(B))$ holds
for any $A,B\in\mathcal M$. Particularly, for any $A,B$ with $AB=0$,
we have $\xi(\delta(A)B+A\delta(B))=\delta(\xi AB)=\delta(0)=0$.
Thus
 $\delta(A)B+A\delta(B)=0$ holds for any $A,B$ with $AB=0$, that is,
 $\delta$ meets the condition for $\xi=0$. Then, by (2) of Claim 5, there
 exists an additive derivation $\varphi$ such that
 $\delta(A)=\varphi(A)+\delta(I)A$ for all $A$. Furthermore, as
 $\delta (\xi I)=\xi \delta (I)I+\xi I\delta (I)=2\xi\delta(I)$, we
 see that
 $\delta(\xi I)\in{\mathcal Z}({\mathcal M})$ by Claim 3 and $\varphi(\xi I)=\xi\delta (I)$.
Since $\delta$ is additive, for any rational real number $r$ and
any $A\in\mathcal M$ we have $\delta (rA)=r\delta(A)$. As
$0=-\varphi(I)=\varphi
(i^2I)=\varphi(iI)iI+iI\varphi(iI)=2i\varphi(iI)$, we see that
$\varphi(iI)=0$, which implies that $\delta(iI)=i\delta(I)$ and
hence $\delta (rI)=r\delta(I)$ holds for any rational complex
number $r$.

Thus, if $\xi$ is a  rational complex number, then
$\xi\delta(I)=\delta(\xi I)=2\xi\delta(I)$, which forces
$\delta(I)=0$. Hence $\delta =\varphi$ is an additive derivation.

{\bf Claim 7.} The statements (1), (2), (3) of the theorem hold.

Note that $L(A)=\delta(A)+AS-SA$ for all $A\in{\mathcal M}$ and
$L(I)=\delta(I)$. Hence, by Claims 5 and 6, $L$ has the forms stated
in the theorem. The proof of the theorem is finished. \hfill$\Box$

By checking the proof of Theorem 3.1, for nonrational complex number
$\xi$, $L$ is an additive derivation under some conditions.

{\bf Corollary 3.3.} {\it Let ${\mathcal M}$ be a von Neumann
algebra without central summands of type $I_1$.   Suppose that
$L:{\mathcal M}\rightarrow{\mathcal M}$ is an additive map and $\xi$
is a nonrational complex number. If $L$ satisfies that
$L([A,B]_\xi)=[L(A),B]_\xi+[A,L(B)]_\xi$ for any $A,B\in{\mathcal
M}$ with $AB=0$, and if $L$ satisfies one of the following
additional conditions, then $L$ is an additive derivation with
$L(\xi I)=0$.}

(i) {\it $L$ is continuous when restricted on ${\mathbb C}I$.}

(ii) {\it There exists a positive number $c$ and a subsequence of
integers $k_n\in {\mathbb Z}$ with $|k_n|\rightarrow\infty$ as
$n\rightarrow\infty$ such that $\|L(\xi^{k_n}I)\|\leq
c|\xi|^{k_n}$.}

{\bf Proof.} By Theorem 3.1 there exists an additive derivation
$\varphi$ such that $\varphi(\xi I)=\xi L(I)$ and
$L(A)=\varphi(A)+L(I)A$ for all $A$. Since
$0=\varphi(i^2I)=2i\varphi(iI)$, we have $\varphi(iI)=0$ and
$L(iI)=iL(I)$. Thus
$L(iA)=\varphi(iA)+iL(I)A=\varphi(iI)A+i\varphi(A)+iL(I)A=iL(A)$ for
any $A\in{\mathcal M}$. It follows from the additivity of $L$ that
$$L(rA)=rL(A) \eqno(3.31)$$
holds for any rational complex number $r$ and any $A\in\mathcal M$.
Note that $L(\xi I)=2\xi L(I)$.

Assume that $L$ meets the condition (i). Take rational complex
numbers $r_n$ so that $\lim_{n\rightarrow\infty}r_n=\xi$. Then
$$
0=\lim_{n\rightarrow\infty}L((r_n-\xi)I)=\lim_{n\rightarrow\infty}(r_n-2\xi)L(I)=-\xi
L(I),$$ which implies that $L(I)=0$. So, $L=\varphi$ is an additive
derivation and $L(\xi I)=0$.

Assume that $L$ satisfies the condition (ii). By induction, it is
easily checked that, for any $k\in{\mathbb Z}$,
$$L(\xi^k I)=(k+1)\xi^kL(I).$$
Thus, by the condition (ii), we have
$$ \|L(I)\|=\frac{1}{|k_n+1||\xi|^{k_n}}\|L(\xi^{k_n} I)\|\leq \frac{c}{|k_n|-1}$$
holds for any integer $k_n$, which implies that $L(I)=0$ as
$|k_n|\rightarrow\infty$. Hence $L$ is an additive derivation.
Moreover, $L(\xi I)=2\xi L(I)=0$.\hfill$\Box$

Finally, let us consider the question of characterizing $\xi$-Lie
derivations for $\xi\not=1$. Obviously, if $\xi=0$, then an additive
$\xi$-Lie derivation is an additive derivation; if $\xi=-1$, then an
additive $\xi$-Lie derivation is an additive Jordan derivation. For
the case that $\xi\not=0,1$, by Theorem 3.1, we have

{\bf Corollary 3.4.} {\it Let ${\mathcal M}$ be a von Neumann
algebra without central summands of type $I_1$.   Suppose that
$L:{\mathcal M}\rightarrow{\mathcal M}$ is an additive map and $\xi$
is a scalar with $\xi\not=0,\pm 1$. Then $L$ is a $\xi$-Lie
derivation if and only if $L$ is an additive derivation.}

{\bf Proof.} The ``if'' part is obvious. For the ``only if'' part,
by Theorem 3.1, $L(I)\in{\mathcal Z}({\mathcal M})$ and there exists
an additive derivation $\varphi$ with $\varphi (\xi I)=\xi L(I)$
such that $L(A)=\varphi(A)+L(I)A$  for all $A\in{\mathcal M}$. Since
$L$ is a $\xi$-Lie derivation, we have
$$L(I)-L(\xi I)=L([I,I]_\xi)=[L(I),I]_\xi+[I,L(I)]_\xi=2L(I)-2\xi L(I).\eqno(3.32)$$
Note that $L(\xi I)=\varphi(\xi I)+\xi L(I)=2\xi L(I).$ So,  by
Eq.(3.32), we see that $L(I)=0$ and $L$ is a derivation.
\hfill$\Box$


\end{document}